\documentclass{article}
\usepackage[utf8]{inputenc}
\usepackage{authblk}
\usepackage[T1]{fontenc}

\usepackage{stmaryrd}
\usepackage{amssymb}
\usepackage{color, amsmath,amssymb, amsfonts, amstext,amsthm, latexsym}

\usepackage{amssymb, epsfig, amssymb, latexsym}

\usepackage{amsmath}
\usepackage{graphicx}
\usepackage{longtable}
\usepackage{graphicx}
\usepackage{subfigure}
\usepackage{color}
 \usepackage{cite}

\renewcommand{\phi}{\varphi}

\newtheorem{example}{Example}

\allowdisplaybreaks[4]

\title{Nonparametric inference of stochastic differential equations based on the relative entropy rate}
\author[a]{Min Dai}
\author[b]{Jinqiao Duan}
\author[c]{Jianyu Hu   \thanks{Corresponding author: jyhu94@outlook.com}}
\author[c]{Xiangjun Wang}
\affil[a]{School of Science, Wuhan University of Technology, Wuhan 430070, China}
\affil[b]{Department of Applied Mathematics, College of Computing, Illinois Institute of Technology, Chicago, Illinois 60616, USA}
\affil[c]{School of Mathematics and Statistics, Huazhong University of Science and Technology, Wuhan 430074, China}

\date{\vspace{-5ex}}

\graphicspath{{Images/}}

\begin{document}

\maketitle

\begin{abstract}
The information detection of complex systems from data is currently undergoing a revolution, driven by the emergence of big data and machine learning methodology. Discovering governing equations and quantifying dynamical properties of complex systems are among central challenges. In this work, we devise a nonparametric approach to learn the relative entropy rate from observations of
stochastic differential equations with different drift functions.
The estimator corresponding to the relative entropy rate then is presented via the Gaussian process kernel theory.
Meanwhile, this approach enables to extract the governing equations.
We illustrate our approach in several examples. Numerical experiments show the proposed approach performs well for rational drift functions, not only polynomial drift functions.

\end{abstract}

\section{Introduction}

Dynamical properties of complex systems such as complexity and sensitivity are ubiquitous, which are significant contents in the research of stochastic dynamical systems \cite{Zhang2020StochasticDA,Pantazis2013ParametricSA,Plyasunov2007EfficientSS}. Less well understood, however, is how to better quantify dynamical properties of stochastic dynamics. The relative entropy rate is an effective tool to characterise the complexity and sensitivity of dynamical properties \cite{dupuis2016path,pantazis2013relative}.
Traditional methods are affected by factors such as unknown governing equations. Owing to the decline cost of data storage and computation, as well as the development of machine learning, data-driven discovery methodologies have made great progress. Combining with data-driven discovery methodologies, it is possible to compute the relative entropy rate from time series data.

There exist many different forms in regard to data-driven methods, such as parametric and nonparametric approaches \cite{chen2018neural,Dai2020DetectingTM,lei2016data,rudy2017data,ruttor2013approximate}. The recent sparse identification of nonlinear dynamics method, which proposed by Brunton and co-workers \cite{brunton2016discovering}, is a scriptures of parametric approaches. Opper \cite{opper2017estimator} employs the parametric techniques to compute the relative entropy rate. Whereas some complex systems are too complicated to model precisely via parametric representation. The nonparametric techniques, for example neural networks \cite{Han2018SolvingHP,Ruthotto2019DeepNN}, Gaussian processes \cite{batz2016variational,Ruttor2013ApproximateGP} thus are applied in data-driven modeling. They outperform the parametric methods for obtaining the data-driven models with accuracy and extensiveness, although the parametric methods provide concrete manifestation of models.

The relative entropy or Kullback-Leibler divergence has recently reemerged of machine learning as a cost function describing the difference between two probability distributions \cite{archambeau2008variational,chen2021solving}. It provides a precise characterization for approximating underlying probability distribution by probability distribution from data. The Kullback-Leibler divergence has been applied to infinite dimensional probability measures of stochastic processes, such as probability inference \cite{apinski2015algorithms} and optimal control problems \cite{kappen2005path,kappen2012optimal,kappen2016adaptive}.
The relative entropy rate \cite{lai2009relative}, which is the long-term average limit of the Kullback-Leibler divergence, is a distance measure of two stochastic processes. It plays an important role in quantifying dynamical properties. Dupuis et.al. \cite{dupuis2016path} discussed the uncertainty and sensitivity bounds of stochastic dynamics from observation data using the relative entropy rate. Pantazis and Katsoulakis \cite{pantazis2013relative} applied the relative entropy rate as a suitable information-theoretic object to analyze the sensitivity of the probability distribution of stochastic processes in regard to perturbations in the parameters of the complex dynamics. Techakesari et.al. \cite{Techakesari2013RelativeER} use the relative entropy rate to design hybrid system filters in the presence of (parameterised) model uncertainty.

In this work, we propose a nonparametric approach to learn the relative entropy rate from observations of stochastic differential equations with different drift functions.
Current efforts are mainly focused on deriving a simple estimator for relative entropy rate via the Gaussian process kernel theory. Meanwhile, this approach enables to extract the governing equations.
Our method requires only one sample data of the stochastic differential equation. It performs well for rational drift functions, not only polynomial drift functions.

The remainder of this paper is structured as follows. In section 2, we introduce the relative entropy rate of stochastic differential equations and the variational formulation. In section 3, we propose a nonparametric technique for the calculation of the relative entropy rate based on one path data of stochastic differential equations. Some numerical examples are presented in section 4, followed in section 5 by conclusion.

\section{The relative entropy rate}
In this work, we consider a stochastic differential equation for the dynamics of a $d$-dimensional diffusion process $X_t\in {\mathbb{R}^d}$ given by
\begin{equation}\label{eq:1}
dX_t=g(X_t)dt+\sigma(X_t){dB_t},
\end{equation}
with initial data $X_0=x_0\in\mathbb{R}^d$, where $g(\cdot)\in{\mathbb{R}^d}$ is the drift function, the diffusion function $\sigma(\cdot)$ is the $d\times{k}$ dimensional matrix and $B_t$ is the standard Brownian motion in $\mathbb{R}^k$.

The generator $\mathcal{A}$ of this stochastic differential equation \eqref{eq:1} is \cite{duan2015introduction}
\begin{equation}\label{eq:2}
\mathcal{A}\phi=g\cdot\nabla \phi + \frac{1}{2}\mathrm{tr}[\sigma\sigma^{T}H(\phi)],
\end{equation}
for all $\phi$ is in Sobolev space $H_0^2(\mathbb{R}^d)$, where $H$ is the Hessian operator with $H=\nabla\nabla^{T}$. Thus, the probability density $p(x,t)$ of the solution process $X_t$ satisfies the Fokker-Planck equation
\begin{equation}\label{eq:3}
p_t=\mathcal{A}^{*}p,
\end{equation}
with initial condition $p(x,0)=\delta(x-x_0)$, where $\mathcal{A}^{*}$ is the adjoint operator of the generator $\mathcal{A}$ in Hilbert space $L^2(\mathbb{R}^d)$, given by
\begin{equation}\label{eq:4}
\mathcal{A}^{*}\phi=-\nabla\cdot(g\phi)+\frac{1}{2}{\mathrm{tr}[H(\sigma\sigma^{T}\phi)]}.
\end{equation}
Here, $H(\sigma\sigma^{T}\phi)$ is interpreted as matrix multiplication of $H=\nabla\nabla^{T}$ and $\sigma\sigma^{T}\phi$ (note that $\phi$ is a scalar function). We call $p$ the stationary probability density of the solution process $X_t$ if it satisfies the stationary Fokker–Planck equation $\mathcal{A}^{*}p=0$.

\subsection{The relative entropy rate for stochastic differential equations}
The relative entropy or Kullback–Leibler divergence between the probability measures $P^g$ and $P^r$ of two solution processes for stochastic differential equation (\ref{eq:1}) with different drifts $g$ and $r$ is defined as
\begin{equation}\label{eq:5}
D_T(P^g,P^r)=E_{P^g}\left[\operatorname{\ln}\frac{P^g}{P^r}\right]=\frac{1}{2}{\int_0^T dt\int{p^g(x,t)||g(x)-r(x)||_{D^{-1}}^2}}dx.
\end{equation}
Here, the diffusion matrix of two solution processes is both $D(x)\doteq\sigma(x)\sigma(x)^T$, $p^g(x,t)$ is the probability density of solution process with drift function $g$ and $||u(x)||_A^2\doteq u(x)\cdot A(x)u(x)$ for some positive definite matrix $A$. Suppose that the stationary probability exists. The probability density $p^g(x,t)$ then converges to the stationary probability density $p^g$ as $t\rightarrow{\infty}$. Hence, we consider the relative entropy rate as follows
\begin{equation}\label{eq:6}
    d(P^g,P^r)=\lim_{T\rightarrow{\infty}}\frac{1}{T}D_T(P^g,P^r)=\frac{1}{2}\int{p^g(x)||g(x)-r(x)||_{D^{-1}}^2}dx.
\end{equation}
In our work, we suppose that the drift $r(x)$ and the diffusion $D(x)$ are known, but the expression form of the drift $g(x)$ is unknown. We, however, want to estimate the relative entropy rate via the observation data of the process $X_t$ with drift $g(x)$ on a large time $T$.

Based on equation \eqref{eq:6}, the estimator of the drift $g(x)$ and the stationary probability density $p^g(x)$ are crucial to the calculation of the relative entropy rate. In order to simplify the estimation problem, we assume that the diffusion $D(x)=\sigma^2{I}$ and drift $g(x)$ satisfies a potential condition such as $g(x)=-\nabla\psi(x)$. Then the stationary probability density fulfills $p^g(x)\varpropto{e^{-\frac{2\psi(x)}{\sigma^2}}}$ and one can estimate the density from observation data.

We will introduce a different way to calculate the relative entropy rate based on the variational formulation and generalised potential condition. To be specific, suppose that the expression of drift $g(x)$ is
\begin{equation}\label{eq:7}
    g(x)=r(x)+D(x)\nabla\psi^*(x).
\end{equation}
Then the relative entropy rate \eqref{eq:6} becomes
\begin{equation}\label{eq:8}
    d(P^g,P^r)=\frac{1}{2}\int{p^g(x)||\nabla\psi^*(x)||_D^2}dx.
\end{equation}
Moreover, the stationary Fokker-Planck equation for drift $g$ is defined as
\begin{equation}\label{eq:9}
    \mathcal{A}_g^*{p^g(x)}=\mathcal{A}_r^*{p^g(x)}-\nabla\cdot(D(x)\nabla\psi^*(x)p^g(x))=0,
\end{equation}
where $p^g(x)$ is the stationary probability density and the operator $\mathcal{A}_r^*$, corresponding to known drift $r(x)$, can express as
\begin{equation}\label{eq:10}
    \mathcal{A}_r^*p^g(x)=-\nabla\cdot(r(x)p^g(x))+\frac{1}{2}\mathrm{tr}[\nabla\nabla^T(D(x)p^g(x))].
\end{equation}

\subsection{Variational formulation}
In this section, we will give a brief introduction about the variational formulation. Here, assume that the specific form of drift $g$ is unknown. We first estimate the drift $g$ from the observation data using the variational formulation for the stationary Fokker-Planck equation \eqref{eq:9} through relative entropy rate. Suppose that the stationary probability density $p^g(x)$ is given and we search for a estimator of the drift $g(x)$ by minimizing the relative entropy rate \eqref{eq:6}. Introducing a Lagrange multiplier function $\psi(x)$, we may derive the drift $g(x)$ from the following Lagrange functional
\begin{equation}\label{eq:11}
\begin{split}
    &\frac{1}{2}\int{p^g(x)||g(x)-r(x)||_{D^{-1}}^2}dx-\int\psi(x)\mathcal{A}_g^*{p^g(x)}dx\\
    =&\frac{1}{2}\int{p^g(x)||g(x)-r(x)||_{D^{-1}}^2}dx-\int\psi(x)\{\mathcal{A}_r^*{p^g(x)}-\nabla\cdot((g(x)-r(x))p^g(x))\}dx.
\end{split}
\end{equation}
The Fokker-Planck operator $\mathcal{A}_g^*$ is defined in \eqref{eq:4} for drift $g(x)$ and the operator $\mathcal{A}_r^*$ is in \eqref{eq:10}. For more details see appendix.

Furthermore, making a variation of Lagrange functional \eqref{eq:11} with respect to $g(x)-r(x)$, one can obtain $g(x)-r(x)=D(x)\nabla\psi(x)$. Inserting this result back into \eqref{eq:11}, the variational representation of relative entropy rate for a unknown potential $\psi$ is
\begin{equation}\label{eq:12}
    \varepsilon_g[\psi]=\int\Bigg\{\frac{1}{2}||\nabla\psi(x)||^2_D+\mathcal{A}_r\psi(x)\Bigg\}p^g(x)dx,
\end{equation}
where the generator $\mathcal{A}_r$ is adjoint operator of $\mathcal{A}_r^{*}$, \eqref{eq:10} which satisfies $\int\psi(x)\mathcal{A}^*p(x)dx=\int p(x)\mathcal{A}\psi(x)dx$. In addition, the explicit expression of generator $\mathcal{A}_r$ is
\begin{equation*}
    \mathcal{A}_r\psi(x)=r(x)\cdot\nabla\psi(x)+\frac{1}{2}\mathrm{tr}[D(x)\nabla\nabla^T\psi(x)].
\end{equation*}

We next introduce the variational bound for the Lagrange functional \eqref{eq:12}
\begin{equation}\label{eq:13}
    -\varepsilon_g[\psi]\leq{\frac{1}{2}\int p^g(x)||\nabla\psi^*(x)||_D^2dx},
\end{equation}
where equality is achieved while $\psi=\psi^*$. We are surprised to discover that, from the perspective of variational bound, the minimisation of variational representation \eqref{eq:12} can give us a estimator for the potential $\psi$ and also compute the relative entropy rate
\begin{equation}\label{eq:14}
\begin{split}
    &\psi^*(x)=\arg \min_{\psi(x)}\varepsilon_g[\psi(x)],\\
    &d(P^g,P^r)=-\varepsilon_g[\psi^*(x)].
\end{split}
\end{equation}
According to the form of Lagrange functional $\varepsilon_g[\psi]$, we all know that applying the observation data which are the ergodic samples of the process with the drift $g$, one can estimate every potential $\psi$ under the stationary probability density condition. The estimator of the relative entropy rate, then, can be obtained.

Subsequently, we give a simple proof for variational bound \eqref{eq:13} as follows\\
\begin{align}
    &{\frac{1}{2}\int p^g(x)||\nabla\psi^*(x)||_D^2dx}+\varepsilon_g[\psi] \notag\\
    =&{\frac{1}{2}\int p^g(x)\Bigg(||\nabla\psi^*(x)||_D^2+||\nabla\psi(x)||_D^2+2\mathcal{A}_r\psi(x)}\Bigg)dx  \notag\\
    =&{\frac{1}{2}\int p^g(x)\Bigg(||\nabla\psi^*(x)-\nabla\psi(x)||_D^2+2\nabla\psi(x)\cdot D(x)\nabla\psi^*(x)+2\mathcal{A}_r\psi(x)}\Bigg)dx  \\
    =&{\frac{1}{2}\int p^g(x)||\nabla\psi^*(x)-\nabla\psi(x)||_D^2 dx}+{\int \psi(x)(\mathcal{A}^*_r p^g(x)-\nabla\cdot [D(x)\nabla\psi^*(x) p^g(x)])}dx  \notag\\
    =& {\frac{1}{2}\int p^g(x)||\nabla\psi^*(x)-\nabla\psi(x)||_D^2 dx}  \notag\\
    \geq& {0}. \notag
\end{align}
The first equation comes from the variational representation \eqref{eq:12} and the second is due to the algebraic formula $(a-b)^2=a^2-2ab+b^2$. Owing to the relationship of generator $\mathcal{A}_r$ and its adjoint operator $\mathcal{A}^*_r$, as well as the usage of integration by parts, the third equation is established. Inserting \eqref{eq:9} into it later, we obtain the second to last equation.

Before moving on to next section, let us discuss another stochastic presentation \eqref{eq:14}.  Consider
\begin{equation}\label{eq:16}
    \varepsilon_g[\psi]=-\lim_{T\rightarrow{\infty}}\frac{1}{T}E_{P_g}\Bigg[\ln{\frac{dP^{r+D\psi}}{dP^{r}}}\Bigg].
\end{equation}
The ${dP^{r+D\psi}}/{dP^{r}}$ represents the Radon-Nykodim derivative of two probability measures. According to \eqref{eq:5} and \eqref{eq:6}, we see
\begin{equation}\label{eq:17}
    \varepsilon_g[\psi]=d(P^g,P^{r+D\psi})-d(P^g,P^r).
\end{equation}
Since $d(P^g,P^{r+D\psi})\geq{0}$, we thus have
\begin{equation}\label{eq:18}
    -\varepsilon_g[\psi]\leq d(P^g,P^r).
\end{equation}
The equality holds in \eqref{eq:18} when $P^g=P^{r+D\psi}$. In next section, we will present a nonparametric estimator for the potential $\psi$ and the relative entropy rate.

\section{Nonparametric inference for relative entropy rate}
In this framework, our goal is to construct an nonparametric estimator to the the relative entropy rate. We however know from the above theory of the potential $\psi$ and stationary probability density $p^g$ is essential to calculating the relative entropy rate. Hence, we first need to a smooth estimator for stationary probability density $p^g$,  and then compute the potential $\psi$ from data using a nonparametric approach via the variational representation \eqref{eq:12}.

There are many density estimation methods, among which the kernel density estimation is commonly used. The drawback of the kernel density estimation is nevertheless that it neglects the temporal ordering of observation data, as well as the time scale between them. Here, we employ the empirical distribution to replace the exact probability $p^g$,
\begin{equation}
    \hat{p}(x)=\frac{1}{n}\sum_{i=1}^{n}\delta(x-x_i).
\end{equation}
The data $x_1,x_2,\cdots,x^n$ are random, ergodic samples of the process with the exact probability $p^g$. On the other hand, one usually build a parametric estimator of the potential $\psi$, which is represented by a finite set of basis functions. However, this representation of parametric method is not sufficient to express many functions, such as rational functions. We therefore work on a nonparametric estimate, which is a more general approach to these functions.

We now demonstrate the nonparametric estimate method. Introduce a penalty term firstly, which is selected as a quadratic form $\frac{1}{2}\sum_k \omega^2_k/\lambda_k$, where $\omega_k$ are the weights corresponding to the basis functions of the parametric representation and $\lambda_k$ are hyper-parameters, to regularize the potential estimator due to a limited amount of observational data. We can also consider the penalty term from the perspective of pseudo-Bayesian. The $\exp\{-\frac{1}{2}\sum_k \omega^2_k/\lambda_k\}$ is viewed as a Gaussian prior distribution over weights $\omega_k$. Assume that the parametric of the variational representation \eqref{eq:12} is $\varepsilon_{\rm{emp}}[\psi_\omega]$. At the same time, the $\exp\{-C\varepsilon_{\rm{emp}}[\psi_\omega]\}$ can be interpreted as a likelihood, where $C$ represents the proportional weight between the penalty
term and the observation data. The potential function $\psi$ can thus be treated as a Gaussian process according to the above interpretation. 
Inspired by the viewpoint of Gaussian process, we will convert the parameter representation form into the kernel function form using the kernel trick, which makes the information for $\psi$ expressed fully. To this end, we define
\begin{equation}
    K(x,x^{'})=\sum_k \lambda_k\phi_k(x)\phi_k(x^{'}),
\end{equation}
where $\lambda_k$ and $\phi_k$ are the orthonormal eigenvalues and eigenfunctions, respectively. This is also the covariance kernel of the prior of a Gaussian process about the potentials $\psi$.

The regularized functional of the potential function $\psi$ can be defined as
\begin{equation}\label{eq:21}
    C\sum_{i=1}^{n}\Bigg\{\frac{1}{2}||\nabla\psi(x_i)||^2_D+\mathcal{A}_r\psi(x_i)\Bigg\}+\frac{1}{2}\int\int\psi(x)K^{-1}(x,x^{'})\psi(x^{'})dxdx^{'},
\end{equation}
via the kernel approach. The $K^{-1}(x,x^{'})$ represents the inverse of the kernel operator. In addition, the penalty term in \eqref{eq:21} can be proved to be equivalent to the reproducing kernel Hilbert space norm of the drifts $\psi$ defined by the kernel $K$. We will further derive a specific representation of the drift function estimator.

The variation of \eqref{eq:21} with respect to $\psi$ yields
\begin{equation*}
    Cn\{\mathcal{A}_r^{*}\hat{p}(x)-\nabla\cdot(D(x)\nabla\psi(x)\hat{p}(x))\}+\int K^{-1}(x,x^{'})\psi(x^{'})dx^{'}=0.
\end{equation*}
Multiplying both sides by the operator $K$, one obtains
\begin{equation}\label{eq:22}
    \psi(x)+C\sum_{j=1}^{n}(\mathcal{A}_g[\psi])_{x^{'}}K(x,x^{'})_{x^{'}=x_j}=0,
\end{equation}
where the generator acts on the kernel function $K$ as
\begin{equation}\label{eq:23}
    (\mathcal{A}_g[\psi])_{x^{'}}K(x,x^{'})_{x^{'}=x_j}=(r(x^{'})+D(x^{'})\nabla\psi(x^{'}))\nabla K(x,x^{'})+\frac{1}{2}\mathrm{tr}[D(x^{'})\nabla\nabla^{T}K(x,x^{'})].
\end{equation}
From above equations \eqref{eq:22} and \eqref{eq:23}, we see that if we know $\nabla\psi(x)$ at all observation data $x=x_i$, the potential function $\psi(x)$ can be calculated for all $x$. The key to this work, however, is to evaluate the gradient of the potential $\psi(x)$. We next compute $\nabla\psi(x)$ at all observation data points via performing the gradient of equation \eqref{eq:22} and setting $x=x_i$. Thus, there are a series of linear equations
\begin{equation}
    \nabla\psi(x_i)+C\sum_{j=1}^{n}(\mathcal{A}_g[\psi])_{x^{'}}\nabla_{x}K(x,x^{'})_{x=x_i,x^{'}=x_j}=0.
\end{equation}
We so far obtain the drift function value $\nabla\psi^{*}(x_i)$ at every point. Further, substituting into these drift function value, we can calculate the relative entropy rate
\begin{equation}\label{eq:}
    d(P^g,P^r)=\frac{1}{2n}{\sum_{i=1}^n{||\nabla\psi^*(x_i)||_D^2}}.
\end{equation}

\section{Numerical Experiments}
We begin our verification of nonparametric estimate approach by describing a few experimental examples. Without special emphasis, we apply the radial basis function kernel
\begin{equation}
    K(x,y)=\exp\Bigg(-\frac{(x-y)^{T}(x-y)}{2l^{2}}\Bigg),
\end{equation}
where the length scale $l$ is the hyper-parameter. In this section, we also show the nonparametric learning results of the drift function in order to further illustrate the advantages of this method.

\begin{example}
Consider a scalar stochastic dynamical system with the polynomial drift term
\begin{equation}\label{eq:26}
dX_{t}=(4X_t-4X_t^3-\beta(X_t^2+2X_t+1))dt+\sigma(X_t){dW_{t}},~~~  X_0=x,
\end{equation}
where $W_t$ is the standard Brownian motion and the diffusion function is $\sigma(x)=1$. The known drift term $r(x)=4x-4x^3$ corresponds to a stochastic double-well system. Another drift function, however, which generates the samples through long-term observation is $g_\beta(x)=4x-4x^3+\nabla\psi(x)$, where $\nabla\psi(x)=-\beta(x^2+2x+1)$ depends on $\beta$. The drifts $g$ and $r$ are equal when $\beta=0$.
\begin{figure}[htp]
\subfigure[$\beta=1$]{
\begin{minipage}[]{0.5 \textwidth}
\centerline{\includegraphics[width=8cm,height=6cm]{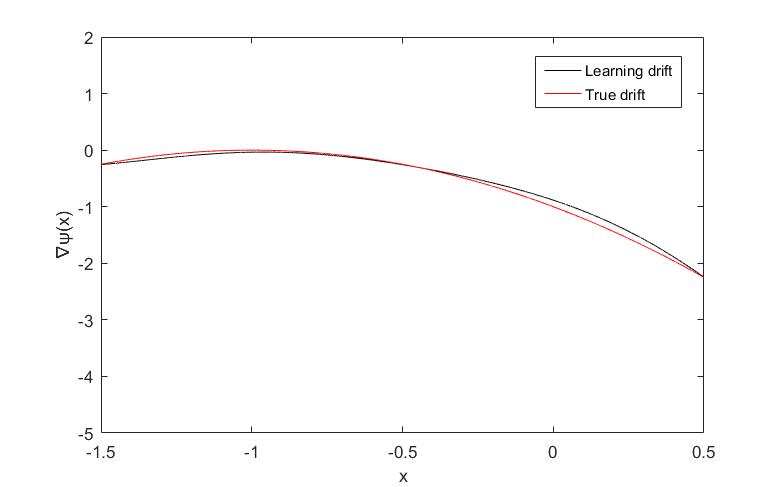}}
\end{minipage}}
\hfill
\subfigure[$\beta=2$]{
\begin{minipage}[]{0.5 \textwidth}
\centerline{\includegraphics[width=8cm,height=6cm]{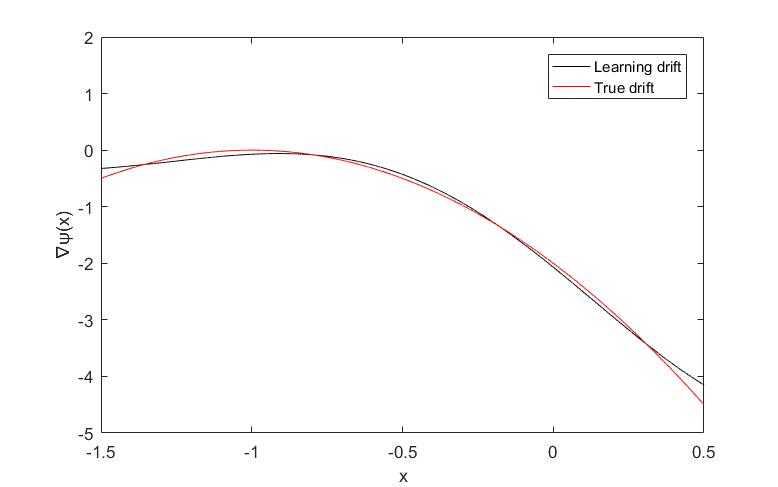}}
\end{minipage}}
\hfill
\subfigure[$\beta=3$]{
\begin{minipage}[]{1.0 \textwidth}
\centerline{\includegraphics[width=8cm,height=6cm]{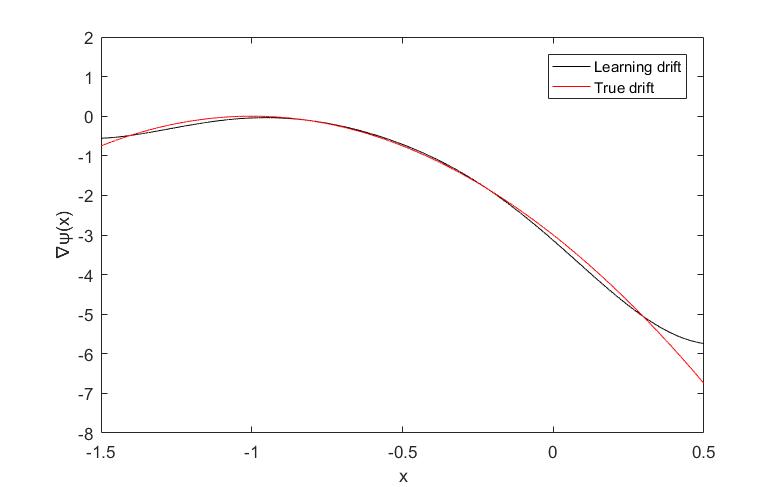}}
\end{minipage}}
\caption{The nonparametric inference results of the drift function $\nabla\psi(x)$. (a) $\beta=1$. (b) $\beta=2$. (c) $\beta=3$.}
\label{Fig1}
\end{figure}

\begin{figure}[htp]
\centering
\includegraphics[width=0.6\textwidth]{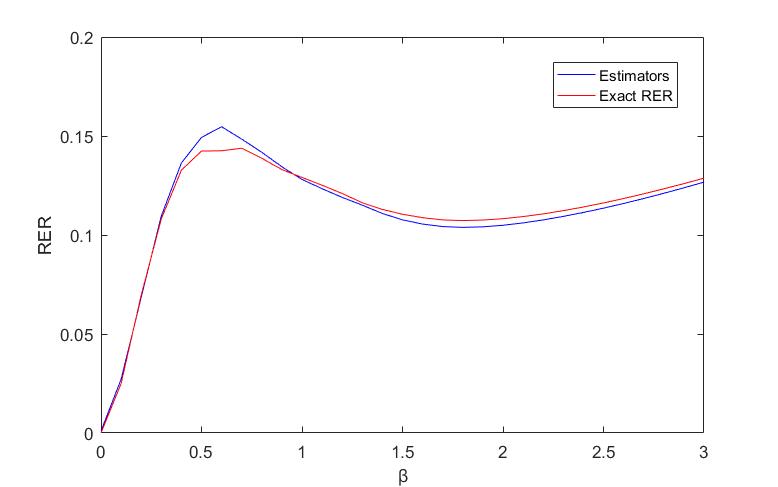}
\caption{The exact relative entropy rate and estimators corresponding to different $\beta$ for the stochastic double-well system.}
\label{Fig2}
\end{figure}
\end{example}
In the numerical simulation, we take the time step $\delta t=0.001$, and use the Euler scheme to sample the stochastic differential equation with drifts $g_\beta$ in regard to the different values of parameter $\beta$. Here, we apply $N=10000$ data points $x_k=X_{t_k}$, uniformly sampled from a trajectory with the observation time length $T=10^5$. Thus the discrete times are $t_k=k\tau$ for $k=1,\cdots,N$ with $\tau=10$. Employing this data points, we can estimate the drifts $\psi(x)$ and their gradients, and further compute the relative entropy rate in term of \eqref{eq:14}. The numerical results will later be presented in Fig.\ref{Fig1} and Fig. \ref{Fig2}.

In Fig.\ref{Fig1}, we compare the accurate and the learning results with regards to potential function $\psi(x)$ for some $\beta$ values and discover that the learning results have a good performance. As shown in Fig.\ref{Fig2}, the exact relative entropy rate and estimators corresponding to different values of parameter $\beta$ are plotted. We can see that the exact relative entropy rate agrees well with the estimate result directly obtained from nonparametric method.

\begin{example}
Consider the transcription factor activator (TF-A) monomer concentration of stochastic differential equation in the gene regulation system \cite{cheng2019most}
\begin{equation}\label{eq:27}
dX_{t}=\beta\Bigg(\frac{6X_t^2}{X_t^2+10}-X_t+0.4\Bigg)dt+\sigma(X_t){dW_{t}},~~~  X_0=x,
\end{equation}
where the diffusion function is $\sigma(x)=1$ and $W_t$ is the standard Brownian motion. The known drift term here is $r(x)=0$, and another drift which generates the samples through long-term observation is $g_\beta(x)=\nabla\psi(x)=\beta((6x^2/(x^2+10))-x+0.4)$ depending on $\beta$. If parameter $\beta=0$, the drifts $g$ and $r$ then are equal.
\begin{figure}[htp]
\subfigure[$\beta=1$]{
\begin{minipage}[]{0.5 \textwidth}
\centerline{\includegraphics[width=8cm,height=6cm]{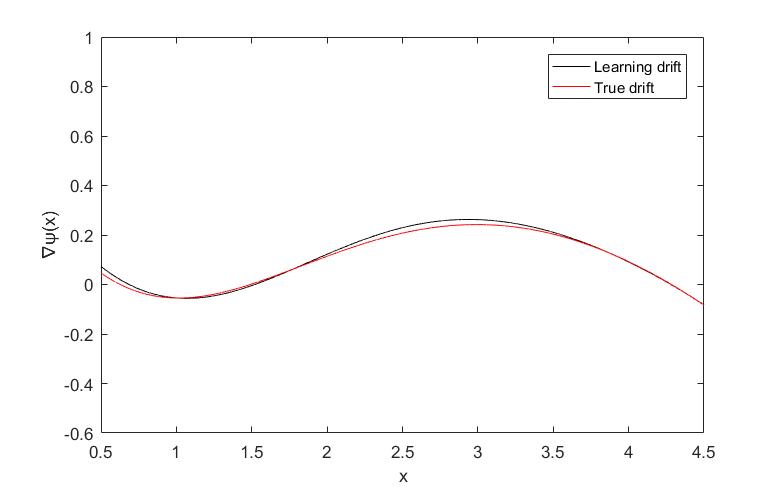}}
\end{minipage}}
\hfill
\subfigure[$\beta=2$]{
\begin{minipage}[]{0.5 \textwidth}
\centerline{\includegraphics[width=8cm,height=6cm]{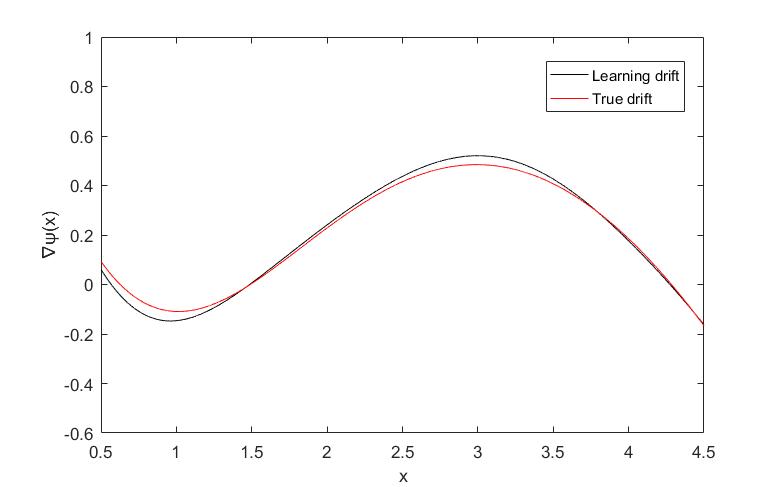}}
\end{minipage}}
\hfill
\subfigure[$\beta=3$]{
\begin{minipage}[]{1.0 \textwidth}
\centerline{\includegraphics[width=8cm,height=6cm]{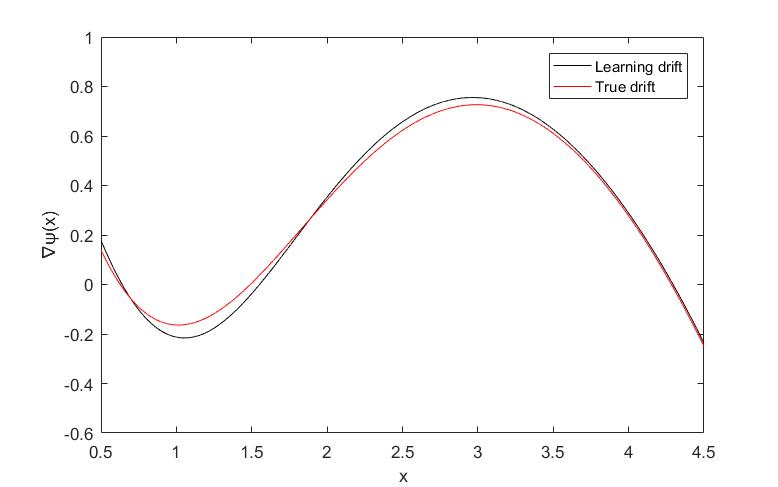}}
\end{minipage}}
\caption{The nonparametric inference results of the drift function $\nabla\psi(x)$. (a) $\beta=1$. (b) $\beta=2$. (c) $\beta=3$.}
\label{Fig3}
\end{figure}

\begin{figure}[htp]
\centering
\includegraphics[width=0.6\textwidth]{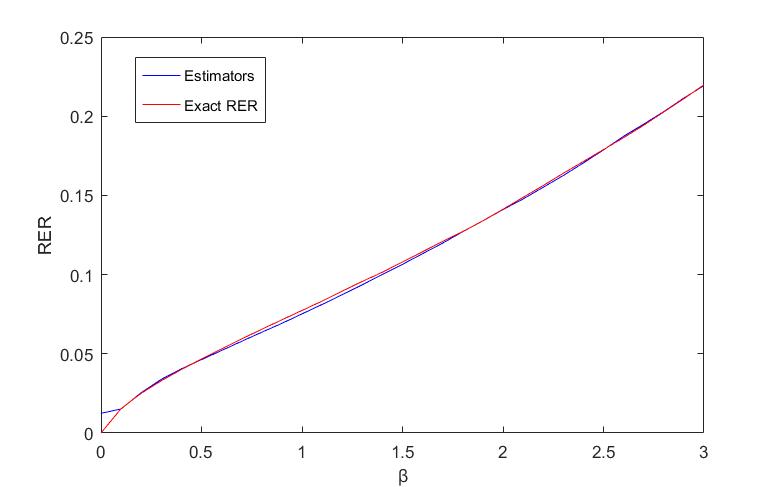}
\caption{The exact relative entropy rate and estimators corresponding to different $\beta$ for the stochastic gene regulation system.}
\label{Fig4}
\end{figure}
\end{example}
To further validate our method, taking the time step $\delta t=0.001$ in the numerical simulation, and using the Euler scheme to sample the stochastic differential equation with drifts $g_\beta$ in regard to the different values of parameter $\beta$, we generate a trajectory with the observation time length $T=10^5$. Here, we utilize $N=10000$ data points $x_k=X_{t_k}$, uniformly sampled from this generated trajectory. Thus the discrete times are $t_k=k\tau, k=1,\cdots,N$ with $\tau=10$. Employing this data points, we can obtain the potentials $\psi(x)$ and their gradients, and meanwhile compute the relative entropy rate in term of \eqref{eq:14}.

As shown in Fig.\ref{Fig3} and Fig.\ref{Fig4}, the comparison regarding the learning and accuracy results of the drift function $\psi(x)$ are given, as well as the exact relative entropy rate and estimators for different values of parameter $\beta$ are plotted. We can see that the evaluation of the drift function $\psi(x)$ is pretty and the exact relative entropy rate agrees well with the estimate result directly obtained from nonparametric method.

\begin{example}
In order to show the superiority of our method, we also consider the trigonometric polynomial function
\begin{equation}\label{eq}
dX_{t}=\beta(sin(X_t)-sin(X_t)^3)dt+\sigma(X_t){dW_{t}},~~~  X_0=x.
\end{equation}
Here, $W_t$ is the standard Brownian motion and the diffusion function $\sigma(x)=1$. Taking the known drift term $r(x)=0$, and another drift which generates the samples through long-term observation $g_\beta(x)=\nabla\psi(x)=\beta(sin(X_t)-sin(X_t).^3)$ depending on $\beta$. If parameter $\beta=0$, the drifts $g$ and $r$ then are equal.
\begin{figure}[htp]
\subfigure[$\beta=1.5$]{
\begin{minipage}[]{0.5 \textwidth}
\centerline{\includegraphics[width=8cm,height=6cm]{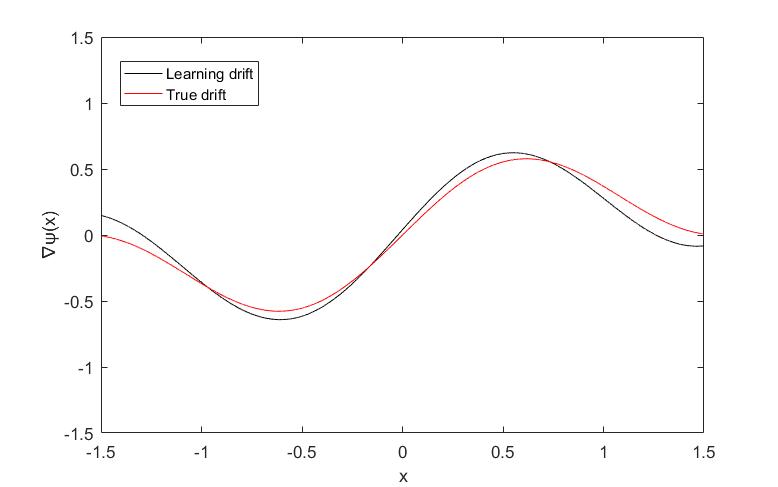}}
\end{minipage}}
\hfill
\subfigure[$\beta=2$]{
\begin{minipage}[]{0.5 \textwidth}
\centerline{\includegraphics[width=8cm,height=6cm]{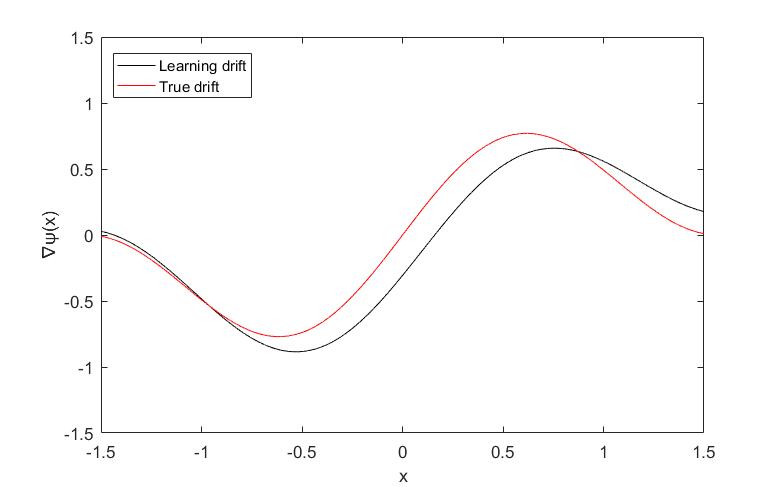}}
\end{minipage}}
\hfill
\subfigure[$\beta=2.5$]{
\begin{minipage}[]{1.0 \textwidth}
\centerline{\includegraphics[width=8cm,height=6cm]{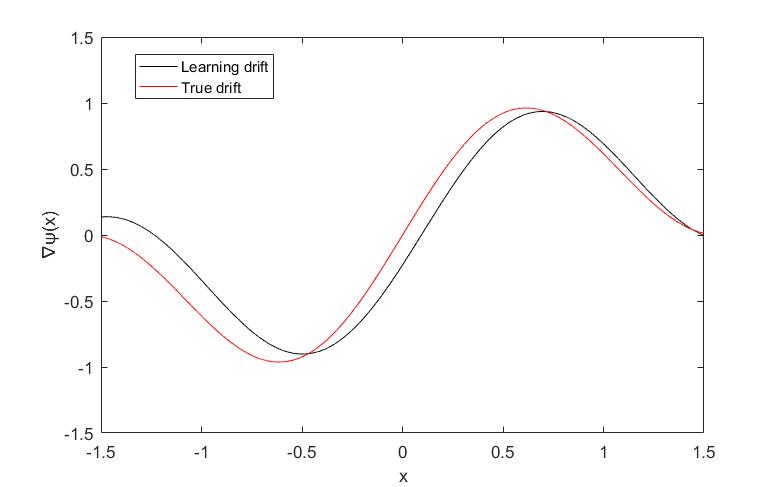}}
\end{minipage}}
\caption{The nonparametric inference results of the drift function $\nabla\psi(x)$. (a) $\beta=1.5$. (b) $\beta=2$. (c) $\beta=2.5$.}
\label{Fig5}
\end{figure}

\begin{figure}[htp]
\centering
\includegraphics[width=0.6\textwidth]{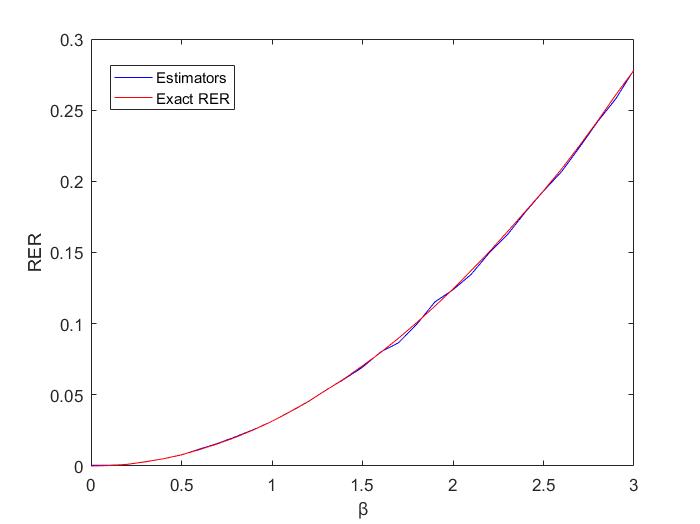}
\caption{The exact relative entropy rate and estimators corresponding to different $\beta$ for the stochastic gene regulation system.}
\label{Fig6}
\end{figure}
\end{example}
We take the time step $\delta t=0.001$ in the numerical simulation, and then apply the Euler scheme to sample the stochastic differential equation with drifts $g_\beta$ in regard to the different values of parameter $\beta$, and generate a trajectory with the observation time length $T=10^5$. Here, we utilize $N=10000$ data points $x_k=X_{t_k}$, uniformly sampled from this generated trajectory, so that the discrete times are $t_k=k\tau, k=1,\cdots,N$ with $\tau=10$. Employing this data points, we obtain the potentials $\psi(x)$ and their gradients, and compute the relative entropy rate in term of \eqref{eq:14}. We show the results in Fig.\ref{Fig5} and Fig.\ref{Fig6}.

The comparison regarding the learning and accuracy results of the drift function $\psi(x)$ are given in Fig.\ref{Fig5}, as well as the exact relative entropy rate and estimators for different values of parameter $\beta$ are plotted in Fig.\ref{Fig6}. We see that the evaluation of the drift function $\psi(x)$ is pretty and the exact relative entropy rate agrees well with the estimate result.

\section{Conclusion and Discussion}
In conclusion, we have presented a nonparametric technique to learn the relative entropy rate which is applied to describe dynamical properties. The study of dynamical properties of systems is significant for stochastic dynamical systems. Here, we offer a good methodology for their study to learn from time series data. Moreover, the nonparametric method not only provides a precise estimation for our functions and the relative entropy rate, but also adapts to a wider range of functions types, such as rational function, compared with parametric approach. This performance has been well demonstrated in examples.

In addition, this work motivates a number of future extensions. Our framework can be extended in other Markov processes, including Markov chains or L\'{e}vy processes. They would be interesting in realistic application. On the other hand, the stationary density is worthy of attention because of the need for a long time observation. We can generalize this issue with finite time window $T$ and evaluate further the marginal densities $p^g(x,t)$.

\section*{Acknowledgments}
This work was supported by the National Natural Science Foundation of China (NSFC) grants 12001213, 11801192, 11771449 and National Science Foundation (NSF) grant 1620449.

\section*{Appendix: Relative entropy}

The relative entropy or Kullback–Leibler divergence between the probability measures $P^g$ and $P^r$ of two solution processes for stochastic differential equation (\ref{eq:1}) with different drifts $g$ and $r$ is defined as
\begin{equation}\label{eq:33}
D_T(P^g,P^r)=E_{P^g}\left[\operatorname{ln}\frac{P^g}{P^r}\right].
\end{equation}
By the Girsanov theorem \cite[Thorem 1.4]{ishikawa2016stochastic} with $\theta=0$, we obtain the Radon–Nykodim derivative
\begin{equation}\label{eq:32}
\frac{dQ}{dP}=\exp\left\{\int_0^Tu(Y_t)dB_t-\frac{1}{2}\int_0^T|u(Y_t)|^2dt\right\}.
\end{equation}
Here, $\sigma(x)u(x)=g(x)-r(x)$, and the stochastic process $Y_t$ is governed by
\begin{equation}\label{eq:34}
dY_t=r(Y_t)dt+\sigma(Y_t){dB_t}.
\end{equation}
Moreover, under the induced probability measure $Q$ by Girsanov theorem, $Y_t$ also satisfies the stochastic differential equation
\begin{equation}\label{eq:36}
dY_t=g(Y_t)dt+\sigma(Y_t){dB^Q_t},
\end{equation}
where $B^Q_t=\int_0^tu(s)ds+B_t$ is a Brownian motion with respect to the induced probability measure $Q$. Then by weak uniqueness of the solutions for stochastic differential equation (\ref{eq:1}) and (\ref{eq:36}), the likelihood function can be represented as
\begin{equation}\label{eq:35}
\operatorname{ln}\frac{P^g}{P^r}=\int_0^Tu(Y_t)dB_t-\frac{1}{2}\int_0^T|u(Y_t)|^2dt.
\end{equation}
To this end, we could use It\^{o} formula to represent the stochastic integral as the a Riemann integral and the relative entropy is
\begin{equation}\label{eq:39}
D_T(P^g,P^r)=\frac{1}{2}\int_0^T\int p^g(t,x)\left[g(x)-{r}(x)\right][D(x)]^{-1}\left[g(x)-{r}(x)\right]dtdx,
\end{equation}
where $D=\sigma\sigma^T$ and $p^g(x,t)$ is the probability density of the stochastic differential equation (\ref{eq:1}) with drift vector filed $g$. Assuming that $p^g(x,t)$ converges to the stationary density $p^g(x)$ for $t\rightarrow\infty$, we shall consider the relative entropy rate \begin{equation}\label{eq:38}
D(P^g,P^r)=\frac{1}{2}\int p^g(x)\left[g(x)-{r}(x)\right][D(x)]^{-1}\left[g(x)-{r}(x)\right]dx.
\end{equation}

Adding the Lagrange multiplier term, equation \eqref{eq:38} becomes our optimal objective functional \eqref{eq:11}.

\bibliographystyle{ieeetr}
\bibliography{ref}

\begin{thebibliography}{10}

\bibitem{Zhang2020StochasticDA}
Z.~guo Zhang, X.~Ma, H.~Yu, and H.~Hua, ``Stochastic dynamics and sensitivity
  analysis of a multistage marine shafting system with uncertainties,'' {\em
  Ocean Engineering}, p.~108388, 2020.

\bibitem{Pantazis2013ParametricSA}
Y.~Pantazis, M.~A. Katsoulakis, and D.~G. Vlachos, ``Parametric sensitivity
  analysis for biochemical reaction networks based on pathwise information
  theory,'' {\em BMC Bioinformatics}, vol.~14, pp.~311--311, 2013.

\bibitem{Plyasunov2007EfficientSS}
S.~Plyasunov and A.~P. Arkin, ``Efficient stochastic sensitivity analysis of
  discrete event systems,'' {\em J. Comput. Phys.}, vol.~221, pp.~724--738,
  2007.

\bibitem{dupuis2016path}
P.~Dupuis, M.~A. Katsoulakis, Y.~Pantazis, and P.~Plech{\'a}c, ``Path-space
  information bounds for uncertainty quantification and sensitivity analysis of
  stochastic dynamics,'' {\em SIAM/ASA Journal on Uncertainty Quantification},
  vol.~4, no.~1, pp.~80--111, 2016.

\bibitem{pantazis2013relative}
Y.~Pantazis and M.~A. Katsoulakis, ``A relative entropy rate method for path
  space sensitivity analysis of stationary complex stochastic dynamics,'' {\em
  The Journal of Chemical Physics}, vol.~138, no.~5, p.~054115, 2013.

\bibitem{chen2018neural}
R.~T. Chen, Y.~Rubanova, J.~Bettencourt, and D.~Duvenaud, ``Neural ordinary
  differential equations,'' {\em arXiv preprint arXiv:1806.07366}, 2018.

\bibitem{Dai2020DetectingTM}
M.~Dai, T.~Gao, Y.~Lu, Y.~Zheng, and J.~Duan, ``Detecting the maximum
  likelihood transition path from data of stochastic dynamical systems.,'' {\em
  Chaos}, vol.~30 11, p.~113124, 2020.

\bibitem{lei2016data}
H.~Lei, N.~A. Baker, and X.~Li, ``Data-driven parameterization of the
  generalized {Langevin} equation,'' {\em Proceedings of the National Academy
  of Sciences}, vol.~113, no.~50, pp.~14183--14188, 2016.

\bibitem{rudy2017data}
S.~H. Rudy, S.~L. Brunton, J.~L. Proctor, and J.~N. Kutz, ``Data-driven
  discovery of partial differential equations,'' {\em Science Advances},
  vol.~3, no.~4, p.~e1602614, 2017.

\bibitem{ruttor2013approximate}
A.~Ruttor, P.~Batz, and M.~Opper, ``Approximate {Gaussian} process inference
  for the drift function in stochastic differential equations,'' in {\em
  Advances in Neural Information Processing Systems}, pp.~2040--2048, Citeseer,
  2013.

\bibitem{brunton2016discovering}
S.~L. Brunton, J.~L. Proctor, and J.~N. Kutz, ``Discovering governing equations
  from data by sparse identification of nonlinear dynamical systems,'' {\em
  Proceedings of the National Academy of Sciences}, vol.~113, no.~15,
  pp.~3932--3937, 2016.

\bibitem{opper2017estimator}
M.~Opper, ``An estimator for the relative entropy rate of path measures for
  stochastic differential equations,'' {\em Journal of Computational Physics},
  vol.~330, pp.~127--133, 2017.

\bibitem{Han2018SolvingHP}
J.~Han, A.~Jentzen, and W.~E, ``Solving high-dimensional partial differential
  equations using deep learning,'' {\em Proceedings of the National Academy of
  Sciences}, vol.~115, pp.~8505 -- 8510, 2018.

\bibitem{Ruthotto2019DeepNN}
L.~Ruthotto and E.~Haber, ``Deep neural networks motivated by partial
  differential equations,'' {\em Journal of Mathematical Imaging and Vision},
  vol.~62, pp.~352--364, 2019.

\bibitem{batz2016variational}
P.~Batz, A.~Ruttor, and M.~Opper, ``Variational estimation of the drift for
  stochastic differential equations from the empirical density,'' {\em Journal
  of Statistical Mechanics: Theory and Experiment}, vol.~2016, no.~8,
  p.~083404, 2016.

\bibitem{Ruttor2013ApproximateGP}
A.~Ruttor, P.~Batz, and M.~Opper, ``Approximate gaussian process inference for
  the drift of stochastic differential equations,'' in {\em NIPS 2013}, 2013.

\bibitem{archambeau2008variational}
C.~Archambeau, M.~Opper, Y.~Shen, D.~Cornford, and J.~Shawe-Taylor,
  ``Variational inference for diffusion processes,'' {\em Advances in Neural
  Information Processing Systems}, vol.~20, pp.~17--24, 2008.

\bibitem{chen2021solving}
X.~Chen, L.~Yang, J.~Duan, and G.~E. Karniadakis, ``{Solving Inverse Stochastic
  Problems from Discrete Particle Observations Using the Fokker--Planck
  Equation and Physics-Informed Neural Networks},'' {\em SIAM Journal on
  Scientific Computing}, vol.~43, no.~3, pp.~B811--B830, 2021.

\bibitem{apinski2015algorithms}
F.~J. Pinski, G.~Simpson, A.~M. Stuart, and H.~Weber, ``Algorithms for
  {Kullback--Leibler} approximation of probability measures in infinite
  dimensions,'' {\em SIAM Journal on Scientific Computing}, vol.~37, no.~6,
  pp.~A2733--A2757, 2015.

\bibitem{kappen2005path}
H.~J. Kappen, ``Path integrals and symmetry breaking for optimal control
  theory,'' {\em Journal of Statistical Mechanics: Theory and Experiment},
  vol.~2005, no.~11, p.~P11011, 2005.

\bibitem{kappen2012optimal}
H.~J. Kappen, V.~G{\'o}mez, and M.~Opper, ``Optimal control as a graphical
  model inference problem,'' {\em Machine learning}, vol.~87, no.~2,
  pp.~159--182, 2012.

\bibitem{kappen2016adaptive}
H.~J. Kappen and H.~C. Ruiz, ``Adaptive importance sampling for control and
  inference,'' {\em Journal of Statistical Physics}, vol.~162, no.~5,
  pp.~1244--1266, 2016.

\bibitem{lai2009relative}
J.~Lai and J.~J. Ford, ``Relative entropy rate based multiple hidden {Markov}
  model approximation,'' {\em IEEE Transactions on Signal Processing}, vol.~58,
  no.~1, pp.~165--174, 2009.

\bibitem{Techakesari2013RelativeER}
O.~Techakesari and J.~J. Ford, ``Relative entropy rate based model selection
  for linear hybrid system filters of uncertain nonlinear systems,'' {\em
  Signal Process.}, vol.~93, pp.~12--22, 2013.

\bibitem{duan2015introduction}
J.~Duan, {\em An Introduction to Stochastic Dynamics}, vol.~51.
\newblock Cambridge University Press, 2015.

\bibitem{cheng2019most}
X.~Cheng, H.~Wang, X.~Wang, J.~Duan, and X.~Li, ``Most probable transition
  pathways and maximal likely trajectories in a genetic regulatory system,''
  {\em Physica A: Statistical Mechanics and Its Applications}, vol.~531,
  p.~121779, 2019.

\bibitem{ishikawa2016stochastic}
Y.~Ishikawa, {\em Stochastic Calculus of Variations}.
\newblock de Gruyter, 2016.

\end{thebibliography}

\end{document}